\documentclass{article}
\usepackage{latexsym}
\usepackage{amssymb}
\usepackage{amsfonts}                     
\usepackage{amsmath}
\usepackage[all]{xy}
\usepackage{color}
\usepackage{mathrsfs}
\usepackage[utf8]{inputenc}
\usepackage[T1]{fontenc}
\usepackage{amscd}
\usepackage{amsxtra} 
\newtheorem{thm}{\bf Theorem}[section]
\newtheorem{cor}[thm]{\bf Corollary}
\newtheorem{lem}[thm]{\bf Lemma}
\newtheorem{prop}[thm]{\bf Proposition}
\newtheorem{defn}[thm]{\bf Definition}

\newcommand{\field}[1]{\mathbb{#1}}

\newcommand{\N }{\field{N}}

\def\WC{{\cal WC}}
\def\SF{{\cal SF}}

\def\SFD{{\rm S\text{-}sfd}}

\def\glsfD{{\rm S\text{-}gl.sf.D}}
\def\glrfD{{\rm gl.sf.D}}

\def\wD{{\rm wD}}

\def\Ext{{\rm Ext}}

\def\proof{{\parindent0pt {\bf Proof.\ }}}

\def\Ker{{\rm ker}}

\def\pd{{\rm pd_R}}
\def\fd{{\rm fd_R}}
\def\id{{\rm id_R}}
\def\glD{{\rm gl.D}}
\def\wD{{\rm w.D}}

\def\SSF{S\text{-}{\cal SF\,}}
\def\SWC{S\text{-}{\cal WC\,}}

\newcommand{\cqfd}
{\hspace{1cm}
\rule{2mm}{2mm}%
\medbreak%
\par%
}
\begin{document}
\title{The strongly flat dimension of  modules and rings}
\author{Ayoub Bouziri}

\date{}
 
\maketitle

\begin{abstract} Let $R$ be a commutative ring with identity, and let $S$ be a multiplicative subset of $R$. Positselski and Slávik introduced the concepts of $S$-strongly flat modules and $S$-weakly cotorsion $R$-modules, and they showed that these concepts are useful in describing flat modules and Enochs cotorsion modules over commutative rings (see the discussion in \cite[Section 1.1]{Pos1}). In this paper, we introduce a homological dimension, called the $S$-strongly flat dimension, for modules and rings. These dimensions measure how far away a module $M$ is from being $S$-strongly flat and how far a ring $R$ is from being $S$-almost semisimple.  The relations between the $S$-strongly flat dimension and other dimensions are discussed. 
\end{abstract}

{\scriptsize \textbf{Mathematics Subject Classification (2020)}: 13D05, 13C15, 13B30.}

 {\scriptsize \textbf{Key Words}: strongly flat modules, strongly flat dimensions, global strongly flat dimensions, weakly Matlis localization.}

\section{Introduction} \hskip .5cm  Throughout this paper, $R$ is a commutative ring with identity, all modules are unitary and $S$ is a multiplicative subset of $R$; that is,  $1 \in S$ and $s_1s_2 \in S$ for any $s_1,s_2 \in S$.  We denote the set of all regular elements and the set of units in $R$ by $\text{reg(R)}$ and $U(R)$, respectively.  Recall that $\text{reg(R)} = \{a \in R : ann_R(a) = 0\}$, where $ann_R(a)$ is the annihilator of an element, and $U(R)$ are examples of multiplicative subsets. Let $M$ be an $R$-module. As usual, we use $\pd(M)$, $\id(M)$ and $\fd(M)$ to denote, respectively, the classical projective dimension, injective dimension and flat dimension of $M$, and, $\wD(R)$ and $\glD(R)$ to denote, respectively, the weak and global homological dimensions of $R$. Also, we denote by $M_S$ the localization of $M$ at $S$. Recall that $M_S \cong M \otimes_R R_S$. 
%In particular, we write $R_S$ to means the ring localization. Furthermore, if $S=\text{reg(R)}$, we denote its total ring of quotients as $Q$.

\medskip         Let $R$ be an integral domain, and let $K$ denote its field of quotients. An $R$-module $C$ is said to be weakly cotorsion if $\Ext^1_R(K, C) = 0$. An $R$-module $F$ is strongly flat if $\Ext^1_R(F, C) = 0$ for all weakly cotorsion $R$-modules $C$. These  concepts first appeared in the paper  \cite{Trl1} and were thoroughly investigated by various researchers (see, for example, \cite{Baz2, Fuc5, lee1}). In the case where $R$ is a commutative ring with zero divisors, weakly cotorsion and strongly flat modules are defined by replacing $K$ with $Q$, the total ring of quotients of $R$.

\medskip
In 2019 Positselski and Slávik introduced the concept of $S$-strongly flat modules and $S$-weakly cotorsion $R$-modules  to describe flat modules and  Enochs cotorsion modules over commutative rings (see the discussion in \cite[Section 1.1]{Pos1}). $S$-weakly cotorsion and $S$-strongly flat $R$-modules  are defined, once again, by replacing $Q$ with $R_S$ in the above definition.  It was shown in \cite[Lemma 1.2]{Baz1} that an $R$-module $F$ is $S$-strongly flat if and only if it is a direct summand of an $R$-module $G$ fitting into a short exact sequence of the type $$ \,\, 0 \to R^{(\alpha)} \to G \to  R_S^{(\beta)} \to 0$$ for some cardinals $\alpha$ and $\beta $. Another interesting result from \cite[Theorem 7.9]{Baz1} states that the class of all $S$-strongly flat $R$-modules is covering if and only if if every flat $R$-module is $S$-strongly flat.
 A general framework for $S$-weakly cotorsion and $S$-strongly flat modules was developed in the paper \cite{Pos1}. We will use $S$-$\SF$ (resp., $S$-$\WC$) to denote the class of $S$-strongly flat modules (resp., $S$-weakly cotorsion modules).

\medskip
%In 2023, Zhang \cite{Zha1} investigated when a cotorsion pair $(\SSF, \SWC)$ is hereditary; that is, $\SSF$ is closed under kernels of monomorphisms, or equivalently, $\SWC$ is closed under cokernels of epimorphisms \cite[Lemma 5.24]{Gob1}. He shows that, for a multiplicative subset consisting of (some) nonzero-divisors in $R$, $(\SSF, \SWC)$ is hereditary if and only if the projective dimension of the $R$-module $R_S$ does not exceed $1$. 

\medskip  
 In this paper, we aim to define a dimension, called the $S$-strongly flat dimension, for modules and rings. It measures how far away a module is from being $S$-strongly flat, and how far away a ring is from being $S$-almost semisimple. Recall from \cite{Bou1} that a ring is said to be $S$-almost semisimple if every $R$-module is $S$-strongly flat.

 \medskip
  The organization of the paper is as follows: In Section $2$ we collect a few known or essentially known results supported by brief arguments and references. We prove that $\SF$ is resolving if and only if the cotorsion pair  $(\SF,\WC)$ is hereditary if and only if the first syzygy module $H$ of the $R$-module $R_S$, i.e., the leftmost term of a short exact sequence of $R$-modules $0\to H \to P \to R_S\to 0$ with a projective $R$-module $P$, is $S$-strongly flat. A multiplicative subset with this property received the name weakly Matlis multiplicative subset. If $S$ consists of (some) nonzero-divisors in $R$, this is equivalent to $\pd(R_S)\leq 1$ \cite[Theorem 3.1]{Zha1}.

In Section 3, we introduce the notion of $S$-strongly flat dimensions and present some general results. For an $R$-module $M$, the $S$-strongly flat dimension $\SFD(M)$ of $M$ is defined as the smallest integer $n \geq 0$ such that $\Ext^{n+1}_R(M, C) = 0$ for any $S$-weakly cotorsion $R$-module $C$. If there is no such $n$, we set $\SFD(M) = \infty$. We then deal with the global $S$-strongly flat dimension, denoted as $\glsfD(R)$, defined as the supremum of the $S$-strongly flat dimensions of $R$-modules. These dimensions exhibit favorable behavior when $\SSF$, the class of all $S$-strongly flat $R$-modules, is resolving. In light of this assumption, several results follow directly from \cite{Son1}. As applications, we provide some characterizations of $S$-almost semisimple rings (see Proposition \ref{3-prop-s-almo-semi}). Furthermore, since $\text{reg(R)}$-almost semisimple rings coincide with semisimple rings, we deduce a new characterization of semisimple rings as well (see Corollary \ref{3-cor-semi}).

\medskip
The concept of $S$-weakly cotorsion dimensions for modules and rings, which can be viewed as the dual of $S$-strongly flat dimensions, paves the way for a forthcoming paper \cite{Bou2}.   This subsequent work presents a distinct set of results on Matlis multiplicative subsets.

\section{Preliminaries}

\hskip .5cm  Let $R$ be a commutative ring, and let $S$ be a multiplicative subset of $R$.

An $R$-module $M$ is said to be $S$-weakly cotorsion if $\Ext^1_R(R_S, M) = 0$, where $R_S$ denotes the localization of the ring $R$ at the multiplicative subset $S$. An $R$-module $F$ is said to be $S$-strongly flat if $\Ext^1_R(F, C) = 0$ for all $S$-weakly cotorsion $R$-modules $C$.  Recall that $(\SSF, \SWC)$ forms a complete cotorsion pair, a fact that follows directly from \cite[Theorem 6.11]{Gob1}.   

\begin{lem}\label{2-lem-close}
\begin{enumerate}
\item  The class of all $S$-weakly cotorsion $R$-modules $S$-$\WC$ is closed under extensions, direct summands and direct products.
\item The class of all $S$-strongly flat $R$-modules $S$-$\SF$ is closed under extensions, direct summands and direct sums.
\end{enumerate}
\end{lem}
%%%%
%%%%%%

\proof Closedness under direct summands is obvious. The assertions about closedness with respect to extensions are given by the long exact sequence. For the direct products and direct sums, we use \cite[Proposition 7.21 and Proposition 7.22]{Rot1}.\cqfd

\begin{lem}[\cite{Baz1}, Lemma 1.2]\label{2-lem-s-strong-flat-charte} An $R$-module $F$ is $S$-strongly flat if and only if $F$ is a direct summand of an $R$-module $G$ fitting into a short exact sequence of the type
\begin{center}
$0\to R^{(\alpha)}\to G\to R_S^{(\beta)}\to 0$
\end{center}
for some cardinals $\alpha$ and $\beta$.
\end{lem}

  \begin{cor}\label{cor-dim-s-strong-flat}
  For any $S$-strongly flat $R$-module $F$, $\pd(F)\leq \pd(R_S).$
  \end{cor}
  
  \proof Follows from Lemma \ref{2-lem-s-strong-flat-charte} and \cite[Exercise 8.6]{Rot1}. \cqfd
%%%%%%%%%%%%%%%%%%%%%%%%%%%%%%%%
%%%%%%%%%%%%%%%%%%%%%%%%%%%%%%%%%%%%%%%

The behavior of strongly flat modules is particularly favorable over Matlis domains. This observation can be traced back to an earlier result by Matlis \cite{Mat1}, which asserts that Matlis domains can be characterized by the property that strongly flat modules form a resolving class; that is, the class of all strongly flat $R$-modules is closed under the kernels of epimorphisms. This was partially the motivation behind the work in \cite{Zha1}, where the focus was on characterizing when $\SSF$ forms a resolving class for a regular multiplicative subset $S$; equivalently, this corresponds to determining when the cotorsion pair $(\SSF, \SWC)$ is hereditary. The key finding was that $\SSF$ is resolving if and only if the projective dimension of $R_S$ is at most $1$. A multiplicative subset with this property is known as a Matlis multiplicative subset. In \cite{Zha1}, the author refers to a ring $R$ as an $S$-Matlis when $S$ is a Matlis multiplicative subset of $R$.

It is essential to note that, in general, the condition $\pd(R_S) \leq 1$ seems to be not necessary for $\SSF$ to be a resolving class. Indeed, consider the first syzygy module $H$ of the $R$-module $R_S$; that is, the leftmost term of a short exact sequence of $R$-modules $0 \to H \to P \to R_S \to 0$, with  a projective $R$-module $P$. Suppose that $H$ is $S$-strongly flat. Then, for all $S$-weakly cotorsion $R$-modules $C$, we have $\Ext^2_R(R_S, C) = \Ext^1_R(H, C) = 0$. Now, let $0\to K\to L\to M\to 0$ be an exact sequence with $K$ and $L$  $S$-weakly cotorsion $R$-modules. The induced exact sequence  
$$\Ext^1_R(R_S,L)=0\to \Ext^1_R(R_S,M)\to  \Ext^2_R(R_S,K)=0 $$  shows that $\Ext^1_R(R_S,M)=0$, and hence, $M$ is $S$-weakly cotorsion. Thus, $\SWC$ is a cosolving class; equivalently, $\SSF$ is a resolving class.  These remarks lead us to introduce the following notion.

\begin{defn}\label{2-def} Let $R$ be a ring and $S$ a multiplicative subset of $R$. Then $S$ is called a weakly Matlis multiplicative subset\footnote{Following the terminology used in \cite{Zha1}, we can also say that a ring $R$ is weakly $S$-Matlis when  $S$ is a weakly Matlis multiplicative subset of $R$.} if $H$, the first syzygy module of the $R$-module $R_S$, is $S$-strongly flat. 
\end{defn}

The following lemma  justify the definition \ref{2-def}:
\begin{lem}\label{2-lem-syzygies}
Let $(K_n)_{n\geq 0}$ and $(K'_n)_{n\geq 0}$ be \textbf{syzygies} of an $R$-module $M$ defined by two projective  resolutions of $M$. Then, for each  $n\in \N$, $K_n$  is  $S$-strongly flat  if and only if $K'_n$ is $S$-strongly flat.
\end{lem}
\proof This follows by the following isomorphism $\Ext^1_R(K_n, N) = \Ext^1_R(K'_n, N)$ where $N$ is an 
 arbitrary $R$-module  \cite[Proposition 8.5]{Rot1}.\cqfd

%Recall that a cotorsion pair $(\X, \Y)$ is said to be hereditary if whenever $0\to X' \to X \to X'' \to 0$ is exact with $X, X''\in \X$, $X'$ is also in $\X$; equivalently, whenever $0\to Y' \to Y \to Y'' \to 0$ is exact with $Y, Y'\in \Y$, $Y''$ is also in $\Y$ \cite[Lemma 5.24]{Gob1}.

 The following lemma can be considered, in some way, as an extension of  \cite[Theorem 3.1]{Zha1}.  
%%%%%%%%%%%%%%%%%%%%%%%%%%%%%%
\begin{lem}\label{2-lem-(sf,wc)-herditary}  Let $R$ be a commutative ring, and $S$ a multiplicative subset of $R$.  The following conditions are equivalent. 
\begin{enumerate}
\item $S$ is a weakly Matlis multiplicative subset.
\item  $\SSF$ is resolving.
\item $\SWC$ is coresolving.
\item The cotorsion pair $(\SSF,\SWC)$ is hereditary. 
\item $\Ext_R^n(F,C)=0$ for any $F\in \SSF$, $C\in \SWC$ and any $n\geq 1$.
\end{enumerate}
\end{lem}
%%%%%%%%%%%%%%%%%
\proof 
$1.\Rightarrow 3.$  See the discussion before the Definition \ref{2-def}. 

 $2.\Rightarrow 1.$  Obvious. 
 
  $2.\Leftrightarrow 3.\Leftrightarrow 4.\Leftrightarrow 5.$ Follows from \cite[Lemma 5.24]{Gob1}
 \cqfd

  \begin{cor}
Assume that $\pd( R_S)\leq 1$. Then $\SSF$ is  resolving. Moreover, if $S$ is consists of (some) nonzero-divisors in $R$, then the converse holds true.
\end{cor}
\proof  The last statement follows from \cite[Corollary 3.2]{Zha1}. \cqfd

%%%%%%%%%%%%%%%%%%%%%%%%%%%%%%%%%
%%%%%%%%%%%%%%%%%%%%%%%%%%%%%%%%%%%Let  $0\to H \to P \to R_S\to 0$  be an exact sequence with a projective $R$-module $P$. Suppose that the $R$-module $H$ is $S$-strongly flat. Then we have $\Ext^2_R(R_S, C) = \Ext^1_R(H, C) = 0$ for all $S$-weakly cotorsion $R$-modules $C$.  Let $0\to K\to L\to  M\to 0$ be an exact sequence with $K$ and $L$ are $S$-weakly cotorsion $R$-module. Then, the induced exact sequence $$\Ext^1_R(R_S,L)=0\to \Ext^1_R(R_S,M)\to  \Ext^2_R(R_S,K)=0 $$ show that $\Ext^1_R(R_S,M)=0$, and hence, $M$ is $S$-weakly cotorsion.
%%%%%%%%%%%%%%%%%%%%%%%%%%%%%%%%%%%%%

\section{$S$-Strongly flat dimension}
\hskip .5cm Let $M$ be an $R$-module. We define the $S$-strongly flat dimension of $M$, denoted by $\SFD(M)$, as the smallest non-negative integer $n$ such that
$$ \Ext^{n+1}_R(M, C) = 0$$ for all $S$-weakly cotorsion $R$-modules $C$.  If there is no such non-negative integer $n$ satisfying the condition, then the dimension is considered to be infinite, and we set $\SFD(M) = \infty$. The $R$-module $M$ is $S$-strongly flat if and only if its $S$-strongly flat dimension is $0$, i.e., $\SFD(M) = 0$.

\begin{thm}\label{3-thm-sfd(m)-char} Let $R$ be a commutative ring, and let $S$ be a multiplicative subset of $R$. Consider the following conditions:
\begin{enumerate}
\item $\SFD(M) \leq n$.
\item $\Ext^{n+1}_R(M, C) = 0$ for any $S$-weakly cotorsion $R$-module $C$.
\item $\Ext^{n+i}_R(M, C) = 0$ for any $S$-weakly cotorsion $R$-module $C$ and any $i \geq 1$.
\item There exists an exact sequence $0 \to F_n\to \cdots\to F_1 \to F_0 \to M\to 0$ with $F_0, F_1, . . . , F_n$ $S$-strongly flat.
\item $F_n$ is $S$-strongly flat whenever there exists an exact sequence $0 \to F_n \to\cdots \to F_1 \to F_0 \to M \to 0$ with $F_0, F_1, . . . , F_{n}$ $S$-strongly flat.
\item Every projective resolution of $M$ has its $n$-th syzygy $S$-strongly flat.
\item There exists a projective resolution of $M$ whose $n$-th syzygy is $S$-strongly flat.
\end{enumerate}

The implications $(3) \Rightarrow (2)\Leftrightarrow  (6)\Leftrightarrow  (7) \Rightarrow (1),$ $(5)\Rightarrow  (7)\Rightarrow (4)$ hold true. Assuming that $S$ is a weakly Matlis  multiplicative subset, all the seven conditions are equivalent.
\end{thm}
\proof The implications $(3) \Rightarrow (2)\Rightarrow (1)$, $(5)\Rightarrow  (7)\Rightarrow (4) $ are evident. 
 
 $(6) \Leftrightarrow (7)$. Follows from Lemma \ref{2-lem-syzygies}. 

$(2) \Leftrightarrow (7)$  Consider an exact sequence: 
$$\xymatrix{ 0 \ar[r]& F_{n} \ar[r]^{} & P_{n-1}   \cdots  \ar[r] &P_{1} \ar[r]^{d_1}& P_{0} \ar[r]^{d_0}& M\ar[r]&0  }$$ be an exact sequence with $P_0, P_1, ..., P_{n-1}$ projective. The equivalence can be deduced from the isomorphism:
$$\Ext^1_R(F_n, C)\cong \Ext^{n+1}_R(M, C),$$ for any $R$-module $C$  \cite[Proposition 8.5]{Rot1}. 

Now, if $S$ is a weakly Matlis multiplicative subset, then, by Lemma \ref{2-lem-(sf,wc)-herditary}, the pair $(\SSF, \SWC)$ forms a hereditary cotorsion pair. Consequently, the equivalences of $(1), (2), (3), (4),$ and $(5)$ can be derived from \cite[Theorem 2.1]{Son1}. \cqfd

We deduce the following inequalities, revealing the connections between projective, $S$-strongly flat, and projective dimension of modules:

\begin{cor}\label{3-cor-eq-fd-sfd-pd} Let $R$ be a commutative ring, and let $M$ be an $R$-module. Then, we have the following:
\begin{enumerate}
\item $\fd(M) \leq   \SFD(M) \leq \pd(M)$
\item $\pd(M) \leq   \SFD(M) + \pd(R_S)$.
\end{enumerate}
\end{cor}
\proof 
$1.$ Assume that $  \SFD(M)=n<\infty$. Then $\Ext_R^{n+1}(M, C)=0$ for any $S$-weakly cotorsion $R$-module $C$. By Theorem \ref{3-thm-sfd(m)-char}, there exists a projective resolution of $M$ whose $n$-th syzygy is $S$-strongly flat. Hence, since every $S$-strongly flat module is flat, $M$ has a flat resolution of length $n$. Hence, $\fd(M)\leq n$. Similarly, one can prove $  \SFD(M) \leq \pd(M)$.

$2.$ We may assume that $  \pd(R_S)=n<\infty$ and $\SFD(M)=m< \infty $. By the proof of the equivalence   $2.\Leftrightarrow 7.$ in Theorem \ref{3-thm-sfd(m)-char}, there exists an exact sequence $$0\to F_m\to P_{m-1}\to P_{m-2} \to \cdots\to P_0\to M\to 0$$ with $P_i$ projective and $F_m$ $S$-strongly flat.  By Corollary \ref{cor-dim-s-strong-flat}, $\pd(F_m)\leq \pd(R_S)=n$. Hence, $\pd(M)\leq m+ n$. \cqfd

\begin{prop} Let $S$ be a weakly Matlis multiplicative subset. Let $0 \to A \to B \to C \to 0$ be an exact sequence of $R$-modules.
If two of  $\SFD(A), \SFD(B)$, $\SFD(C)$ are finite, so is the third. Moreover,
\begin{enumerate}
\item  $\SFD(B) \leq \sup\{\SFD(A), \SFD(C)\}$.
\item  $\SFD(A) \leq \sup\{\SFD(B), \SFD(C)-1\}$.
\item  $\SFD(C) \leq \sup\{\SFD(B), \SFD(A)+1\}$.
\item $\SFD(C) = \SFD(A) +1$ in case $B$ is $S$-strongly flat and $C$ is not.
\end{enumerate}
	\end{prop}

\proof By Lemma \ref{2-lem-(sf,wc)-herditary}, this is a special case of \cite[Proposition 2.4]{Son1}. \cqfd

As is customary for global dimensions, we define the strongly flat dimensions of commutative ring $R$ with respect to a multiplicative subset $S$ of $R$:

\begin{defn}
Let $R$  be a commutative ring and $S$ be a multiplicative subset of $R$. The global $S$-strongly flat dimension of $R$, denoted by $\glsfD(R)$, is defined as the supremum of the $S$-strongly flat dimensions of all $R$-modules. When $S = \text{reg}(R)$, we simply write $\glrfD(R)$.
\end{defn}
Notice that if $S$ is a multiplicative subset of $R$ such that $R_S$ is projective as $R$-module (for example, $S = U(R)$), then $\glsfD(R) = \glD(R)$.

 Recall from \cite[Definition 7.6]{Baz1} that a commutative ring $R$ is called $S$-almost perfect  if $R_S$ is a perfect ring and $R/sR$ is a perfect ring for every $s \in S$. According to \cite[Theorem 7.9]{Baz1}, a commutative ring $R$ is  $S$-almost perfect if and only if every flat $R$-module is $S$-strongly flat. Thus, if $R$ is $S$-almost perfect, then the global $S$-strongly flat dimension and the weak global dimension coincide. Moreover,
\begin{prop}
Let $R$ be a Noetherian commutative ring. Assume that $R$ is either of finite spectrum or of Krull dimension $1$. Then there exists a Matlis multiplicative subset $S$ of $R$ such that for any $R$-module $M$, $\fd(M)=\SFD(M)$. Consequently, $\wD(R) = \glsfD(R)$. 
\end{prop}

\proof In the case where $R$ has a finite spectrum, there exists a multiplicative subset $S$ of the form $\{s, s^2, \ldots, s^n, \ldots\}$ such that flat modules and strongly flat modules coincide \cite[Theorem 1.12]{Pos1}. Moreover, $\pd(R_S)\leq 1$ by \cite[Lemma 1.9]{Pos3}.

For the case where $R$ has Krull dimension $1$, the multiplicative subset is chosen to be $S=R\setminus\bigcup\limits_{q\in P_0}q$, where $P_0$ is the set of all minimal prime ideals \cite[Theorem 13.9]{Pos2}. Moreover, $\pd(R_S)\leq 1$ by \cite[Corollary 13.7]{Pos2}.\cqfd

\begin{thm}\label{3-main-thm} Let $R$ be a commutative ring and let $S$ be a multiplicative subset of $R$. The the following assertions hold true:
\begin{enumerate}
\item $\wD(R)\leq \glsfD(R)= \glD(R)$

\item  $\glD\leq \glsfD(R)+\pd(R_S)$.
% Moreover, if $\pd(Q)=1$ then $\glD=  \glsfD(R)+1$
 
\item If $S$ is weakly Matlis, then  $\glsfD(R) \\
                = \sup\{\SFD(M)| \text{$M$ is a finitely generated $R$-module}\}\\
                = \sup\{\SFD(M)| \text{ $M$ is a cyclic $R$-module}\}\\
                = \sup\{\id(C)|\text{ C $S$-weakly cotorsion}\}\\
                = \sup\{\SFD(C)|\text{ C $S$-weakly cotorsion}\}.$

\item If  $S$ is weakly Matlis  and $\glsfD(R) < \infty$, then  $\glsfD(R)\\ = \sup\{\id(M) |\text{ $M$ is a $S$-weakly cotorsion $S$-strongly flat $R$-module}\}\\
= \sup\{\SFD(M)| \text{ $M$ is an injective $R$-module}\}.$
\end{enumerate}
\end{thm}

\proof $(1)$ and $(2)$ follow from Corollary \ref{3-cor-eq-fd-sfd-pd}.

By Lemma \ref{2-lem-(sf,wc)-herditary} $(\SSF,\SWC)$ is a hereditary cotorsion pair. Hence,  $(3)$ and $(4)$ follow from \cite[Theorem 3.1]{Son1}. \cqfd

%\begin{thm} Let $R$ be a Noetherian commutative ring with finite spectrum.  \end{thm}

\begin{cor}\label{3-cor-n}  Let $S$ be a weakly Matlis  multiplicative subset. Then the following are equivalent for an integer $n \geq 0$:
\begin{enumerate}

\item $\glsfD(R) \leq n$.
\item Every cyclic $R$-module is of $S$-strongly flat dimension  $\leq n$.
\item All $S$-weakly cotorsion $R$-modules are of injective dimension $\leq n$.
\item All $S$-weakly cotorsion $R$-modules are of $S$-strongly flat dimension $\leq n$.
\item $\glsfD(R) \leq \infty$, and all $S$-strongly flat $S$-weakly cotorsion $R$-modules are of injective dimension $\leq n$.
\item $\glsfD(R) \leq \infty$, and all injective $R$-modules are of $S$-strongly flat dimension $\leq n$.
\item $\Ext^{n+1}_R(M, N) = 0$ for all $S$-weakly cotorsion $R$-modules $M$ and $N$.
\item $\Ext^{n+j}_R(M, N) = 0$ for all $S$-weakly cotorsion $R$-modules $M$, $N$ and $j\geq 1$.
\end{enumerate}
\end{cor}

%\proof Follows from Theorem \ref{3-main-thm}.  \cqfd

It is well known that $\glD(R)\leq1$ if and only if every submodule of a projective $R$-module is projective. Similarly, we demonstrate in the following result that $\glsfD(R)\leq 1$ if and only if every submodule of any $S$-strongly flat $R$-module is $S$-strongly flat.

\begin{cor}  Assume that $S$ is a weakly Matlis  multiplicative subset. Then the following assertions are equivalent:
\begin{enumerate}
\item $\glsfD(R)\leq 1$.
\item Every ideal is $S$-strongly flat.
\item All $S$-weakly cotorsion $R$-modules are of injective dimension $\leq 1$.
\item All $S$-weakly cotorsion $R$-modules are of $S$-strongly flat dimension $\leq 1$.
\item Every submodule of projective $R$-module is $S$-strongly flat.
\item Every submodule of $S$-strongly flat $R$-module is $S$-strongly flat.
\end{enumerate}
\end{cor}
\proof  The equivalences of $(1), (2), (3)$ and $(4)$ follow by setting $n=1$ in Corollary \ref{3-cor-n}. Moreover, the application  $(6)\Rightarrow (5)$ and $(5)\Rightarrow (1)$ are evident.   To complete the proof, it suffices to prove the implication $(1)\Rightarrow (6)$.  Suppose $(1)$ holds. Let $F$ be an $S$-strongly flat $R$-module, and  let $K$ be a submodule of $F$. Consider the following exact sequence $$0\to K\to F\to F/K\to 0.$$   Since, by $(1)$, $\SFD(F/K)\leq 1$, $K$ is $S$-strongly flat by Theorem \ref{3-thm-sfd(m)-char}.  \cqfd

Whenever a kind of global dimension of rings is studied, it is of special interest to characterize a ring $R$ with such dimension zero. As far as $\glsfD(R)$ is concerned, we have the following result which can be seen as a special case of \cite[Theorem 3.3]{Son1}.
% when $R$ is assumed to be weakly $S$-Matlis.  

Recall from \cite{Bou1} that a commutative ring $R$ is said to be $S$-almost semisimple if every $R$-module is $S$-strongly flat. 
\begin{prop}\label{3-prop-s-almo-semi} Let $R$ be a commutative ring and let $S$ be a multiplicative subset of $R$. Then the following assertions  are equivalent.
\begin{enumerate}
\item $R$ is $S$-almost semisimple.
\item Every $R$-module is $S$-strongly flat (i.e, $\glsfD(R) = 0$).
\item Every $S$-weakly cotorsion $R$-module is injective.
\item $\pd(R_S)\leq 1$ and every $S$-weakly cotorsion $R$-module is $S$-strongly flat.
\item $\pd(R_S)\leq 1$ and $\Ext^1_R(M, N) = 0$ for all $S$-weakly cotorsion $R$-modules $M$ and $N$.
\item $\pd(R_S)\leq 1$ and $\Ext^n_R(M, N) = 0$ for all $S$-weakly cotorsion $R$-modules $M$, $N$ and $n \geq 1$.
\item  Every $S$-weakly cotorsion $R$-module has an injective envelope with the unique mapping property.
\item $\pd(R_S)\leq 1$ and every $S$-weakly cotorsion $R$-module has an $S$-strongly flat cover with the unique mapping property.
\end{enumerate}
\end{prop}

\proof We only verify the proof.  The implications $(1)\Leftrightarrow (2)\Leftrightarrow (3)$, $(4)\Leftrightarrow (5)$, $(4)\Rightarrow (8)$, $(3)\Rightarrow (7)$ 
 are clear. 

$(1)\Leftrightarrow (4)$ Follows by \cite[Theorem 2.7]{Bou1}.
 
$(5)\Leftrightarrow (6)$ Follows by Lemma \ref{2-lem-(sf,wc)-herditary}.

$(7)\Rightarrow (3)$  Let $C$ be an $S$-weakly cotorsion $R$-module. Consider an exact sequence 

$$\xymatrix{ 0 \ar[r]& C \ar[r]^{f} & E_1 \ar[r]^{g}   &E_{2} },$$ where $E_1$ and $E_2$ are injective and $f$ is an injective envelope of $C$ with the unique mapping property. We have $gf = 0 = 0 f$. It follows that $g= 0$. Then $C$ is isomorphic to $E_1$  and therefore it is injective. 

$(8)\Rightarrow (4)$  Let $C$ be an $S$-weakly cotorsion $R$-module. By (8), there exists  an exact sequence 

$$\xymatrix{  F_2 \ar[r]^{g} & F_1 \ar[r]^{f}   & C \ar[r]^{}&0 },$$ where $F_1$ and $F_2$ are $S$-strongly flat and $f$ (resp., g) is an  $\SSF$-cover of $C$ (resp., $\Ker(f)$) with the unique mapping property (notice that $\ker(f)$ is $S$-weakly cotorsion by Wakamutsu's Lemmas \cite[Lemma 2.1.1]{Xu1}). Since $fg=0=f0$, $g=0$, and hence, $ker(f)=0$. Therefore, $C\cong F_1$ is $S$-strongly flat. \cqfd

%\begin{rem} The equivalences of $(2)$, $(6)$ and $(7)$ of Proposition \ref{3-prop-s-almo-semi} appeared in \cite[Theorem 2.8]{Bou1} without the hypothesis that $R$ is $S$-almost semisimple. \end{rem}

Since, by \cite[Corollary 2.3]{Bou1}, $\text{reg(R)}$-almost semisimple rings are nothing but semisimple rings, we have:

\begin{cor}\label{3-cor-semi}  The following are equivalent for a commutative ring $R$. 
\begin{enumerate}
\item $R$ is semisimple.
\item Every $R$-module is strongly flat (i.e, $\text{gl.sf.D}(R) = 0$).
\item Every weakly cotorsion $R$-module is injective.
\item $R$ is a Matlis ring and every weakly cotorsion $R$-module is strongly flat.
\item $R$ is a Matlis ring and $\Ext^1_R(M, N) = 0$ for all weakly cotorsion $R$-modules $M$ and $N$.
\item $R$ is a Matlis ring and $\Ext^n_R(M, N) = 0$ for all weakly cotorsion $R$-modules $M$, $N$ and $n \geq 1$.
\item  Every weakly cotorsion $R$-module has an injective envelope with the unique mapping property.
\item $R$ is  a Matlis ring and every weakly cotorsion $R$-module has a strongly flat cover with the unique mapping property.
\end{enumerate}
\end{cor}

Ayoub Bouziri: Faculty of Sciences, Mohammed V University in Rabat, Rabat, Morocco.

\noindent e-mail address: ayoub$\_$bouziri@um5.ac.ma; ayoubbouziri25@gmail.com

\end{document}